\documentclass[11pt]{article}
\usepackage{mac}
\renewcommand{\IN}{\mathbf{N}}

\begin{document}
\title{{\Large{\bf A Sub-Gaussian Berry-Esseen Theorem For the 
\\
Hypergeometric
Distribution}}}
\maketitle
\begin{center}
S.N. Lahiri, A. Chatterjee, T. Maiti\\
Department of Statistics\\
Iowa State University\\
Ames, IA 50011\\[1in]
{\bf ABSTRACT}\\
\end{center}
In this paper, we derive a necessary and sufficient condition on the
parameters of the Hypergeometric distribution for weak convergence to a
Normal limit. We establish a Berry-Esseen theorem for the
Hypergeometric distribution solely under this necessary and sufficient
condition. We further derive a nonuniform Berry-Esseen bound where the
tails of the difference between the Hypergeometric and the 
Normal distribution functions are shown
to decay at a  sub-Gaussian rate.\\[.8in]
{\it AMS(2000) Subject Classification} Primary 60F05;   
 Secondary 60G10, 62E20, 62D05.\\
Research partially supported by NSF grant no. DMS 0306574.
\newpage
\section{Introduction}
Consider a dichotomous finite population of size $N$ having $M$
individuals of type A and $N-M$ individuals of type B. Suppose
 a sample
of size $n$ is drawn at random, without replacement from this
population. Let $X$ denote the number of `type A'-individuals 
in the sample. Then, $X$ is said to have  the Hypergeometric 
distribution with parameters
$n,M,N$,  written as $X\sim Hyp(n;M,N)$. The probability mass function
(p.m.f) of $X$ is given by,
\begin{equation}
P(X=x)\equiv P(x;n,M,N)=\left\{\begin{array}{ll}
\frac{
{M\choose x}{N-M\choose n-x}}{{N\choose n}}
\quad &\mbox{if}\quad x=0,1\ldots,n\\ 
0\quad &
\mbox{otherwise,}\end{array}\right.
\label{hyppmf}
\end{equation}
where, for any two integers $r\geq 1$ and $s$, 
\begin{eqnarray}
{r\choose s}=\left\{\begin{array}{ll}
\frac{r!}{s!(r-s)!}\quad
&\mbox{if}\quad 0\leq s\leq r
\\ 0\quad &
\mbox{otherwise,}\end{array}\right.
\label{factorialdef}
\end{eqnarray}
with  $0!=1$ and $r!=1\cdot2\cdot\cdot\cdot r$.
Let $f=\frac{n}{N}$ denote the sampling fraction
and let  $p=\frac{M}{N}$ denote the  proportion
of the  `type A'-objects in the population.
The Hypergeoemetric distribution plays an important role
in many areas of Statistics, including sample surveys
(e.g., finite population inference), statistical quality control
(acceptance sampling plans), etc.  Normal approximations 
to the Hypergeometric probabilities $P(.;n,M,N)$ of
(\ref{hyppmf}) are classical in the cases where the sampling fraction $f$
and the proportion $p$ are bounded away from 0 and 1; see for example
Feller(1971). However, the extreme cases where  $f$ or
$p$ take values near the boundary values $0$ and $1$ are 
very important in sample surveys and quality control 
applications.
In this paper, we investigate
the validity and the rate of Normal approximation to the 
Hypergeometric distribution allowing the parameters 
$f$ and $p$ to tend  to any points in  
the interval $[0,1]$, including  the boundary points. 
The main results of the paper
 give a necessary and sufficient condition on the parameters 
$f$ and $p$  for a valid Normal approximation.
 It is shown that a Normal limit for properly centered
and scaled version of $X$ holds if and only if
\begin{equation}
Np(1-p)f(1-f)\raw \infty.
\end{equation}
As a consequence of this, we  conclude that 
 for the Normal distribution function  to approximate 
the distribution
function of $X$, all four quantities, namely, 
(i)  the number  $M$ ($=Np$) of `type A'-objects,
(ii) the number of `type B'-objects, $N-M$, 
(iii) the sample size $n$,
as well as (iv) the size of the unselected objects $N-n$
in the population,
 must tend to infinity.

We also investigate the  rate of Normal approximation 
to the distribution of $X$. Note that $X$ is the sum
of a collection of $n$ {\it dependent} Bernoulli random 
variables. In Section 2, 
we establish a Berry-Esseen Theorem on  the rate of
 Normal approximation to the distribution function
 of $X$ solely under the necessary and sufficient  
condition  (1.3).  
It is shown that under (1.3)
 the rate of approximation is $O([Np(1-p)f(1-f)]^{-1/2})$.
It is also shown in Section 2 that 
this  rate  is optimal and can not be
 improved.
Note that the rate 
$O([Np(1-p)f(1-f)]^{-1/2})$
 is equivalent to the standard rate  $O(n^{-1/2})$
(for sums of $n$ {\it independent}  Bernoulli random 
variables, say) only when   $p$ is bounded away from 0  and 1 
and $f$ bounded away from $1$. However, for $p$ and
$f$ close to these boundary points, the rate of approximation
can be substantially slower. In such situations, 
the dependence of the Bernoulli random 
variables associated with $X$  has a nontrivial effect 
on the accuracy of the Normal approximation.

Under somewhat  stronger conditions on $f$ and $p$,
  we also derive  a non-uniform version of the Berry-Esseen 
Theorem. The nonuniform bound shows that in the tails, the   
error of Normal approximation  dies at a sub-Gaussian rate
 for a wide range of values of $f$ and $p$. 
As a corollary,  we also derive  an exponential
(sub-Gaussian) probability inequality for the tails of $X$,
which  may be of independent interest.


The rest of the paper is organized as follows. We conclude Section 1
with a brief literature review. Section 2 
introduces the asymptotic framework and
contains the results on the
validity of the Normal approximation and  the Berry-Esseen
theorems. Proofs of all the  results
are given in Section 3. 

For results on 
Normal approximations to Hypergeometric probabilities
in the standard  cases where the sampling fraction $f$
and the proportion $p$ are bounded away from 0 and 1, see 
Feller(1971).  For general $p$ and $f$, 
Nicholson (1956) derived some very precise bounds
for the point probabilities $P(.;n,M,N)$ using some 
nonstandard normalizations of the
Hypergeometric random variable $X$. 
General methods for proving the CLT for sample means under
 sampling  without
replacement from finite populations are 
given by Madow (1948), Erdos \& Renyi
(1959) and Hajek(1960). For results on
Berry-Esseen Theorems and Edgeworth expansions for the functions of sample
means and U-statistics based on finite population observations,
see  Babu \& Singh (1985), Kokic \& Weber (1990), Chen \& Sitter
(1993), Bloznelis (1999), 
Bloznelis \& G\"{o}tze (2000), and the references therein.

\section{Main Results}
\setcounter{section}{2}
\setcounter{equation}{0}
Let $r$ be a positive integer valued variable and for each $r\in
\IN$ (where $\IN=\{1,2,\ldots\}$), let $X_r$ be a random variable 
having the
Hypergeometric distribution with parameters $(n_r,M_r,N_r)$. Thus we
consider a sequence of dichotomous finite populations indexed by $r$,
with the population of objects of type A and the sampling fraction
respectively given by, 
\begin{eqnarray}
p_r=\frac{M_r}{N_r}\qmq{and} f_r=\frac{n_r}{N_r}
\forall r\in \mathbf{N}.
\label{p&fdef}
\end{eqnarray}
To avoid trivialities, all through the paper, we shall assume that
\begin{eqnarray}
1\leq M_r <N_r,\quad 1\leq n_r <N_r \quad \forall r 
\in \mathbf{N},
\qmq{and}
N_r^{-1} = o\left(1\right)\quad r\rightarrow \infty.
\label{ineq}
\end{eqnarray}
Thus, $p_r$ , $f_r \in (0,1)$ for all $r \in \mathbf{N}$.
Let
\begin{equation}
\sigma^2_r \equiv N_rp_rq_rf_r(1-f_r),
\label{avari}
\end{equation}
 where $q_r=1-p_r$. The
first result concerns the validity of the Normal approximation to the
distribution of $X_r$.\\[.2in]
{\bf Theorem 2.1:}~ Suppose that (\ref{ineq}) holds and that $X_r\sim
Hyp(n_r,M_r,N_r)$, $r\in \mathbf{N}$. Then there 
exists a Normal random variable
$W\sim N(\mu,\sigma^2)$ for some $\mu\in \mathbf{R}$ and $\sigma \in
(0,\infty)$ such that 
\begin{eqnarray}
\Delta_r\equiv 
\sup_{x \in \mathbf{R}} \bigg|P\left(\frac{X_r-n_rp_r}{\sigma_r}\leq
x\right)-P\left(W\leq x\right)\bigg|\longrightarrow 0
\qmq{as}r\raw\infty
\label{deltadef}
\end{eqnarray}
if and only if 
\begin{eqnarray}
\sigma^2_r\rightarrow \infty\qquad \textrm{as} 
\quad r\rightarrow \infty.
\label{rcondn}
\end{eqnarray}
When (\ref{rcondn}) holds, one must have $\mu=0$ and $\sigma=1$.\\[.2in]

Note that $\sigma^2_r=n_rp_rq_r(1-f_r)=\frac{N_r-1}{N_r}Var(X_r)$. Hence
Theorem 2.1 shows that the Normal approximation to the 
Hypergeometric distribution holds solely under the condition that the
variance of the Hypergeometric distribution goes to infinity with
$r$. In particular, it is not necessary to impose separate conditions on
the asymptotic behavior of the three sequences $\{n_r\}_{\{r\geq 1\}}$,
$\{p_r\}_{\{r\geq 1\}}$ and $\{f_r\}_{\{r\geq 1\}}$. 
A necessary condition for
(\ref{rcondn}) is that $n_r\rightarrow \infty$ and $(N_r-n_r)\rightarrow
\infty$ as $r\rightarrow \infty$. This follows by noting that
$\sigma^2_r=n_rp_rq_r(1-f_r)=(N_r-n_r)p_rq_rf_r
\leq \min\{n_r,N_r-n_r\}$ for all $r\geq 1$. Thus, for the
Normal approximation to hold, both the sample size $n_r$ and the
residual sample size $(N_r-n_r)$ must become unbounded
as $r\raw\infty$. 
By
interchanging the roles of $p_r$ and $q_r$ with $f_r$ and $(1-f_r)$, it
follows that for the validity of the 
 Normal approximation, we must also have
\begin{eqnarray}
M_r \wedge (N_r-M_r) \longrightarrow \infty\quad \textrm{as} \quad
r\rightarrow \infty,
\label{Scondn1}
\end{eqnarray}
i.e., 
 the number of objects of
type A and type B must go to infinity with $r$. 

In a seminal paper, Hajek
(1968) obtained a necessary and sufficient condition 
for the CLT for finite population sums, assuming that 
\begin{eqnarray}
n_r \wedge  N_r-n_r\rightarrow \infty \quad\textrm{as}
 \quad r\rightarrow \infty.
\label{Hcondn}
\end{eqnarray}
The observations above imply that this is not 
a serious restriction; Indeed, in the cases where 
(\ref{Hcondn}) fail, the CLT need not hold.

 Condition (\ref{rcondn}) also allows the proportion 
$p_r$ of `type A'-objects in the population and the 
sampling fraction $f_r$ to 
simultaneously converge to
the extreme points 0 and 1 at certain rates. If the sequence
$\{f_r\}_{\{r\geq 1\}}$ is bounded away from 0 and 1
and (2.2) holds, then
the CLT of Theorem 2.1 
holds if and only if (iff)
\begin{eqnarray}
\frac{1}{N_r}=o (q_r\wedge p_r)\qmq{as}r\raw\infty,
\label{Scondn2}
\end{eqnarray}
i.e., iff (\ref{Scondn1}) holds. Similarly, for
$\{p_r\}_{\{r\geq 1\}}$ bounded away from 0 and 1, the CLT holds
{\it{iff}} 
\begin{eqnarray}
\frac{1}{N_r}= o( f_r\wedge (1-f_r))\qmq{as}r\raw\infty,
\label{Scondn3}
\end{eqnarray}
i.e., iff  (\ref{Hcondn}) holds. However, when both
$\{p_r\}_{\{r\geq 1\}}$ and $\{f_r\}_{\{r\geq 1\}}$ simultaneously
converge to some limits in $\{0,1\}$, neither of (\ref{Scondn2}) and 
(\ref{Scondn3}) alone is enough to guarantee the CLT. For example if
$f_r \sim N^{-a}_r$ and $p_r \sim N^{-b}_r$ for some $0<a,b<1$, with
$a+b>1$, then (\ref{Scondn2}) and (\ref{Scondn3}) hold but the Normal
approximation of Theorem 2.1 is no longer valid.

Next we obtain a refinement of (\ref{deltadef}) by specifying the rate
of convergence of $\Delta_r$ to zero.\\[.2in]
{\bf Theorem 2.2:}~
Suppose that $X_r\sim Hyp(n_r,M_r,N_r)$, $r\in \mathbf{N}$, and that
(\ref{rcondn}) holds. Then there exists a constant
$C_1\in (0,\infty)$ such that for all $r\in\IN$, 
\begin{eqnarray}
\Delta_r\leq \frac{C_1}{\sigma_r}.
\label{bdelta1}
\end{eqnarray}
\\[.2in]

Theorem 2.2 is a uniform Berry-Esseen theorem that shows that
under (\ref{rcondn}), the  rate of Normal approximation 
to the Hypergeometric distribution is uniformly 
$O\left(\sigma^{-1}_r\right)$ as $r\rightarrow
\infty$. When both the sequences $\{p_r\}_{\{r\geq 1\}}$ and
$\{f_r\}_{\{r\geq 1\}}$ are bounded away from 0 and 1, this rate is
$O\left(n_r^{-\frac{1}{2}}\right)$, which is the same as the rate of
Normal approximation for sums of $n_r$ {\it independent and identically
distributed} (iid) random variables with a finite third moment. 
Although the Hypergeometric random
variable $X_r$ can be written as a sum of $n_r$ {\it dependent}
Bernoulli ($p_r$) variables, the lack of independence 
of the  summands does
not affect the  rate of Normal approximation as long as 
 the sequence $\{p_r\}_{r\geq 1}$ is bounded away from
$0$  and $1$ and  $\{f_r\}_{r\geq 1}$ is bounded away from
$1$; The rate becomes worse otherwise.

A second important aspect of Theorem 2.2  is that the bound on
$\Delta_r$ holds under the same condition (\ref{rcondn}) that is both
necessary and sufficient for a Normal limit. Since
$\frac{X_r-n_rp_r}{\sigma_r}$ is supported on a lattice with maximal
span $\sigma^{-1}_r$, it is not difficult to show that
 if (\ref{rcondn})
holds, then 
$\liminf_{r \rightarrow \infty} \Delta_r\sigma_r >0
$, i.e.,
there exists a constant $C_2 \in (0,\infty)$ such that 
\begin{eqnarray}
\Delta_r>\frac{C_2}{\sigma_r}
\label{bdelta2}
\end{eqnarray}
for all but finitely many $r$'s.
 Thus, the rate in Theorem 2.2
is {\it optimal} and  can not be improved upon. 


The next result gives a non-uniform version of the
Berry-Esseen theorem.
To state it, let $\phi(\cdot)$ and $\Phi(\cdot)$
 respectively denote the density and the distribution function 
of a standard Normal random variable, i.e., 
$
\phi(x)=\frac{1}{\sqrt{2\pi}}\exp(-\frac{x^2}{2}),
\quad x \in \mathbf{R}$ and 
$\Phi(x)=\int_{-\infty}^{x}\phi(t)dt,\quad x \in \mathbf{R}
$.
Also let $I(\cdot)$ denote the indicator function.
Define 
\begin{equation}
\delta_r=\frac{1}{10}{(\mbox{max}(a_{1r},2))}^{-1},\qquad r\geq 1,
\end{equation}
where $a_{1r}=\frac{\bar{f}_r+4}{4(1-\bar{f}_r)}$ and where 
\begin{displaymath}
\bar{f}_r=\left\{\begin{array}{c@{\quad:\quad}l} f_r 
&\textrm{if}\quad f_r \leq \frac{1}{2}   \\
1-f_r & \textrm{if}\quad f_r > \frac{1}{2}.\\\end{array}\right.
\end{displaymath} 
Then, we have the following result.
\\[.2in]
{\bf Theorem 2.3:}~
Suppose that 
$X_r\sim Hyp(n_r,M_r,N_r), r
\in \mathbf{N}$. Assume that $r$ is such that 
\begin{eqnarray} 
\delta_r\sigma_r > 1.
\label{bdelta3}
\end{eqnarray}
Then there exists universal constants $C_3,C_4 \in (0,\infty)$
(not depending  on $r,n_r,M_r$ and $N_r$)  such that 
\begin{eqnarray}
\bigg|P\left(\frac{X_r-n_rp_r}{\sigma_r}\leq x\right)-\Phi(x)\bigg|\leq
\frac{C_3}{\sigma_r}\frac{1
+{|x|}^2}{\lambda_r(x)}\exp\left(-C_4x^2\lambda^2_r(x)\right)
\label{main3bound}
\end{eqnarray}
for all $x \in \mathbf{R}$, where
 $\lambda_r(x)=q_rI(x\leq 0)+p_rI(x\geq 0)$.
\\[.2in]

Theorem 2.3 shows that the error of Normal approximation to the
Hypergeometric distribution dies at a sub-Gaussian rate in
the tails. The only condition needed for the validity of this
bound is (\ref{bdelta3}).  It is easy to
check that
\begin{eqnarray}
\delta_r \in \left(\frac{1}{25},\frac{1}{20}\right]
\label{deltainterval}
\end{eqnarray}
for all $r$ satisfying (\ref{bdelta3}). 
Hence, the bound in (\ref{main3bound}) is available 
for all $r$ such that
$\sigma_r\geq 25$.

An immediate consequence of Theorem 2.3
is the following exponential (sub-Gaussian) probability bound on the
tails of $X_r$.\\[.2in]
{\bf Corollary 2.4:}~
Suppose that $X_r\sim Hyp(n_r,M_r,N_r),r \in \mathbf{N}$.
Then, there exist universal constants $C_5,C_6 \in
(0,\infty)$ (not depending on $r,n_r,M_r,N_r$) such that for all $r$
satisfying (\ref{bdelta3}),
\begin{eqnarray}
P\left(\bigg|\frac{X_r-n_rp_r}{\sigma_r}\bigg|\geq x\right)\leq
\frac{C_5}{{(p_r\wedge q_r)}^3}\exp\left(-C_6x^2{[p_r\wedge
q_r]}^2\right) \forall x\in (0,\infty).\nonumber
\end{eqnarray}


\section{Proofs}
\setcounter{section}{3}
\setcounter{equation}{0}

We now introduce some notation and notational convention
 to be used in this section.
For real numbers $x,y$, let $x\wedge y=min\{x,y\}$ and
 $x\vee y=max\{x,y\}$.  Let $\lfloor x\rfloor$ denote 
the largest integer not exceeding $x$, $x
\in \mathbf{R}$.  
For $a \in (0,\infty)$, write $\phi_a(x)=\frac{1}{a}\phi(\frac{x}{a})$
and $\Phi_a(x)=\Phi(\frac{x}{a})$, $x \in \mathbf{R}$, for the density
and distribution functions of a  $N(0,a^2)$ variable.
Write $\phi_a =\phi$ and  $\Phi_a =\Phi$ for $a=1$.
Let 
\begin{equation}
\Delta^*_r(x)=P\left(\frac{X_r-n_rp_r}{\sigma_r}\leq
x\right)-\Phi(x),\quad x\in \mathbf{R}.
\end{equation}
Let $\mathbf{N}=\{1,2,\ldots\}$,
$\mathbf{Z}_+=\{0,1,\ldots\}$ and 
$\mathbf{Z}=\{\ldots,-1,0,1,\ldots\}$.

For notational simplicity, we shall drop the suffix $r$ from notation,
except when it is important to highlight the dependence on $r$. 
Thus, we write $n,M,N$ for $n_r,M_r,N_r$ respectively and set
$p=\frac{M}{N}$, $q=1-p$ and $f=\frac{n}{N}$. 
We shall use $C$ to denote a generic positive constant that does not
depend on $r$. Unless otherwise stated, limits in order symbols are
taken by letting $r\rightarrow \infty$.

For proving the result, we shall frequently make use
of Stirling's approximation (cf. Feller(1971))
\begin{eqnarray}
m!&=&\sqrt{2\pi}e^{-m+\epsilon_m}m^{m+\frac{1}{2}}
 \forall m \in \mathbf{N},
\label{stirling}
\end{eqnarray}
where the error term $\epsilon_m$ admits the bound 
\begin{displaymath}
\frac{1}{12m+1}\leq \epsilon_m \leq \frac{1}{12m}\qquad \forall
 m \in \mathbf{N}.
\end{displaymath}
Also note that for $g(y)=\log{y}$, $y \in (0,\infty)$,
the $k$th derivative of $g$ is given by 
$g^{(k)}(y)=
\frac{{(-1)}^{k-1}(k-1)!}{y^k}$,  $y \in
(0,\infty), \,k\in \mathbf{N}$. Hence, for any $k\in\mathbf{N}$
 and $\delta \in (0,1)$,
\begin{eqnarray}
\Big|g^{(k)}\left(1+x\right)\Big|
	\leq \frac{(k-1)!}{{(1-\delta)}^k}\qquad
\forall\quad 0 \leq |x|<\delta.
\label{derineq}
\end{eqnarray}

For Lemma 3.1, let $X\sim Hyp(n;M,N)$ for a {\it given}
set of  integers $n,M,N
\in \mathbf{N}$ with $1\leq n\leq (N-1)$, $1\leq M\leq (N-1)$.
 Note that this notation is
consistent with our convention of dropping the suffix $r$; $X,n,M,N$ in
Lemma 3.1 would subsequently 
represent $X_r,n_r,M_r,N_r$ for a {\it{fixed}} $r\in
\mathbf{N}$ for which (2.2) holds. Let 
\begin{equation}
x_{k,n}=\frac{x-np}{\sqrt{npq}}\quad \textrm{and}
\quad a_{k,n}=\frac{x_{k,n}}{(1-f)\sqrt{npq}},\quad 0\leq k\leq n,
\label{defxkn}
\end{equation}
 where 
$f=\frac{n}{N}$, $p=\frac{M}{N}$ and  $q=1-p$.
Lemma 3.1 gives a basic
approximation to Hypergeometric probabilities solely under
condition 
(\ref{future}) stated below.
\\[.2in]

{\bf Lemma 3.1}~
Suppose that  $X\sim Hyp(n;M,N)$ for a {\it given}
set of  integers $n,M,N
\in \mathbf{N}$ 
such that 
\begin{eqnarray}
0<f<1,\quad 0<p<1 \quad\textrm{and} \quad 6(np\wedge nq)\geq 1,
\label{future}
\end{eqnarray}
where  $f=\frac{n}{N}$, $p=\frac{M}{N}$ and  $q=1-p$
are as in (\ref{defxkn}). 
Then, for any given $\delta \in (0,\frac{1}{2}]$, 
\begin{eqnarray}
\log{P(k;n,M,N)}=-\frac{x^2_{k,n}}{2(1-f)}-\frac{1}{2}\log{(2\pi
npq(1-f))}+r^*_n(k)
\label{logexpand}
\end{eqnarray}
for all  $k\in\{0,\ldots,n\}$ 
with $|a_{k,n}|\leq \delta$, where $P(k;n,M,N)=P(X=k)$ (cf. 
(\ref{hyppmf})) and where the remainder term $r^*_n(k)$ admits the
bound
\begin{eqnarray}
|r^*_n(k)|&\leq &\frac{1}{6npq(1-\delta)(1-f)}
+\left[\frac{1}{2}|a_{k,n}|+a^2_{k,n}\left\{\frac{1}{4}
+\frac{2\delta}{{(1-\delta)}^3}\right\}\right]\nonumber\\
	& & \quad {}+{|a_{k,n}|}^3npq\left(\frac{f}{4}
+1\right)\left\{\frac{1}{2}+\frac{2(1+\delta)}{{(1-\delta)}^3}\right\},
\label{rbound}
\end{eqnarray}
provided  $|a_{k,n}|\leq \delta$.
\\[.1in]

\noindent
{\bf Proof:}~ 
%
For
 $k \in \{0,1,\ldots,n\}$,
\begin{eqnarray}
P(k,n,M,N)
& = & 
	\frac{{Np \choose k}{Nq \choose n-k}}{{N \choose n}}\nonumber\\
 & = & {n \choose k}p^kq^{n-k}\frac{{\prod\limits^{k-1}_{j=1}(1
	-\frac{j}{Np})}{\prod\limits^{n-k-1}_{j=1}(1
	-\frac{j}{Nq})}}{\prod\limits^{n-1}_{j=1}(1-\frac{j}{N})}\nonumber\\
	 & = &{n \choose k}p^kq^{n-k}\ R(k,n,M,N),\qmq{say.}
\label{pmfbreak}
\end{eqnarray}
First  consider the denominator of $R(k;n,M,N)$. By (3.2), 
\begin{eqnarray*}
{\prod\limits^{n-1}_{j=1}(1-\frac{j}{N})}&=&\frac{N!}{{(N-n)!}{N^n}}\\
		&=&
\frac{{e^{(-N+\epsilon_N)}}{N^{N+\frac{1}{2}}}}{{e^{(-(N-n)
+\epsilon_{N-n})}}{(N-n)^{N-n+\frac{1}{2}}}}\frac{1}{N^n}\\
		&=&\frac{{e^{(\epsilon_N
-\epsilon_{N-n})}}{e^{-n}}}{(1-f)^{N(1-f)+\frac{1}{2}}}.
\end{eqnarray*}
Similarly,  the numerator of $R(k;n,M,N)$ is given
by
\begin{eqnarray*}
{\prod\limits^{k-1}_{j=1}(1
-\frac{j}{Np})}{\prod\limits^{n-k-1}_{j=1}(1-\frac{j}{Nq})}
&=&
{\frac{{e^{-k}e^{\epsilon_{Np}-\epsilon_{Np-k}}}}{{(1
	-\frac{k}{Np})^{Np-k
		+\frac{1}{2}}}}{\frac{{e^{-(n-k)}e^{\epsilon_{Nq}
		-\epsilon_{Nq-n+k}}}}{{(1-\frac{n-k}{Nq})^{Nq-n+k
			+\frac{1}{2}}}}}}\\
& = &
{\frac{{e^{-n}}{e^{\epsilon_{Np}-\epsilon_{Np-k}+\epsilon_{Nq}
	-\epsilon_{Nq-n+k}}}}{{(1-\frac{k}{Np})^{Np-k
	+\frac{1}{2}}}{(1-\frac{n-k}{Nq})^{Nq-n+k+\frac{1}{2}}}}}.
\end{eqnarray*}
Note that by (3.4), 
\begin{eqnarray}
\frac{k}{Np}=f+x_{k,n}\sqrt{\frac{fq}{Np}}
\qmq{and} \frac{n-k}{Nq}=f-x_{k,n}\sqrt{\frac{fp}{Nq}}\, .
\label{fracdefs}
\end{eqnarray}
Hence $R(k;n,M,N)$ can be expressed as
\begin{eqnarray*}
R(k;n,M,N)=
\exp(\epsilon_{Np}-\epsilon_{Np-k}+\epsilon_{Nq}
	-\epsilon_{Nq-n+k}+\epsilon_{N-n}-\epsilon_{N})
	{(1-f)^{N(1-f)+\frac{1}{2}}}\\
	\quad\times\left\{\left(1
		-f-x_{k,n}\sqrt{\frac{fq}{Np}}\right)^{Np\left(1
	-f-x_{k,n}\sqrt{\frac{fq}{Np}}\right)+\frac{1}{2}}\right\}\\
	\quad\times\left\{\left(1-f
		+x_{k,n}\sqrt{\frac{fp}{Nq}}\right)^{Nq\left(1
		-f+x_{k,n}\sqrt{\frac{fp}{Nq}}\right)+\frac{1}{2}}
			\right\}\, .
\end{eqnarray*}
Next write 
\begin{eqnarray}
z_{k,n}&=&\frac{x_{k,n}\sqrt{\frac{fp}{Nq}}}{1-f},\quad
y_{k,n}=\frac{x_{k,n}\sqrt{\frac{fq}{Np}}}{1-f}
\qmq{and}\nonumber\\
\epsilon^*&=&\epsilon_{Np}-\epsilon_{Np-k}
+\epsilon_{Nq}-\epsilon_{Nq-n+k}+\epsilon_{N-n}-\epsilon_{N}.
\label{eps*}
\end{eqnarray}
Then it follows that
\begin{eqnarray}
\log R(k;n,M,N)
&=&
\epsilon^*-\frac{\log(1-f)}{2}-
	\Big(Np(1-f)(1-y_{k,n})
	+\frac{1}{2}\Big)\log(1-y_{k,n})\nonumber\\
& &
\quad {}-\Big(Nq(1-f)(1+z_{k,n})+\frac{1}{2}\Big)\log(1+z_{k,n})\nonumber\\
&\equiv& 
\epsilon^*-\frac{\log(1-f)}{2}-A_1-A_2,\qquad \mbox{say.}
\label{takelog}
\end{eqnarray}
Fix $\delta \in (0,1/2)$. By Taylor's expansion and (\ref{derineq}), 
\begin{eqnarray}
A_1&=&\left(Np(1-f)(1-y_{k,n})+\frac{1}{2}\right)\log(1-y_{k,n})\nonumber\\
	&=&
\left(Np(1-f)(1-y_{k,n})+\frac{1}{2}\right)\left(-y_{k,n}
-\frac{y^2_{k,n}}{2}+r_{1n}(k)\right)\nonumber\\
 &=&
-y_{k,n}\left(Np(1-f)+\frac{1}{2}\right)
-\frac{y^2_{k,n}}{2}\left(\frac{1}{2}-Np(1-f)\right)+r_{2n}(k),
\label{expanA_1}
\end{eqnarray}
where $r_{1n}(k)$ and $r_{2n}(k)$ are remainder terms, defined by the
equality of the successive expressions. By (3.3),
for all $n,k$ satisfying $|y_{k,n}|\leq \delta$,
\begin{eqnarray}
|r_{1n}(k)|& \leq &
\frac{2}{{(1-\delta)}^3}\frac{{|y_{k,n}|}^3}{3!}\quad 
\mbox{and},\nonumber\\
|r_{2n}(k)|& \leq &
\frac{Np}{2}(1-f){|y_{k,n}|}^3+\bigg|Np(1-f)(1-y_{k,n})
+\frac{1}{2}\bigg|\cdot |r_{1n}(k)| .
\label{errbound1}
\end{eqnarray}
By similar arguments,
\begin{eqnarray}
A_2&=&\left[Nq(1-f)(1+z_{k,n})+\frac{1}{2}\right]\log
(1+z_{k,n})\nonumber\\
	&=&
\left(Nq(1-f)+\frac{1}{2}\right)z_{k,n}
	+\frac{z^2_{k,n}}{2}\left[Nq(1-f)-\frac{1}{2}\right]+r_{3n}(k),
\label{expanA_2}
\end{eqnarray}
where  for all $n,k$, satisfying $|z_{k,n}|\leq \delta$, 
\begin{eqnarray}
|r_{3n}(k)|\leq
Nq(1-f)\frac{{|z_{k,n}|}^3}{2}+\bigg|Nq(1-f)(1+z_{k,n})
+\frac{1}{2}\bigg|\cdot
\frac{{|z_{k,n}|}^3}{3{(1-\delta)}^3}\,\, .
\label{errbound2}
\end{eqnarray}
From, (\ref{takelog}),(\ref{expanA_1}) and (\ref{expanA_2}), we have
\begin{eqnarray}
\log R(k;n,M,N)&=& \epsilon^*-\frac{1}{2}\log
(1-f)-\left[\frac{1}{2}\left(z_{k,n}-y_{k,n}\right)
+\frac{z^2_{k,n}}{2}\left\{Nq(1-f)-\frac{1}{2}\right\}\right.\nonumber\\
	& & \quad \left. {}+\frac{y^2_{k,n}}{2}\left\{Np(1-f)
-\frac{1}{2}\right\}+r_{2n}(k)+r_{3n}(k)\right]\nonumber\\
	&=& \epsilon^*-\frac{1}{2}\log
(1-f)-\frac{x^2_{k,n}f}{2(1-f)}+r_{4n}(k)
\label{logexpann}
\end{eqnarray}
where for all $n,k$ satisfying $(|y_{k,n}|\vee |z_{k,n}|)\leq
\delta$,
\begin{eqnarray*}
|r_{4n}(k)|\leq
|r_{2n}(k)|+|r_{3n}(k)|+\frac{1}{2}|y_{k,n}-z_{k,n}|
+\frac{1}{4}\left(y^2_{k,n}+z^2_{k,n}\right) .
\end{eqnarray*}
Next using Stirling's formula on the binomial term, we have
\begin{eqnarray}
\log \left\{{n\choose k}p^kq^{n-k}\right\}
&=&
 \log\left\{\frac{e^{(\epsilon_n-\epsilon_k
			-\epsilon_{n-k})}}{\sqrt{2\pi npq}}\right\}
	-	\left(nq-x_{k,n}\sqrt{npq}+\frac{1}{2}\right)\log
\left\{1-x_{k,n}\sqrt\frac{p}{nq}\right\}\nonumber\\
& &
 \quad {}-
		\left(np+x_{k,n}\sqrt{npq}+\frac{1}{2}\right)\log
\left\{1+x_{k,n}\sqrt\frac{q}{np}\right\}
\nonumber\\
	&\equiv & \epsilon^{**}-\log \sqrt{2\pi npq}-A_3-A_4,\quad
\mbox{say},
\label{binomexpan}
\end{eqnarray}
where $\epsilon^{**}=\epsilon_n-\epsilon_k-\epsilon_{n-k}$. Next write
$\tilde{y}_{k,n}=x_{k,n}\sqrt{\frac{p}{nq}}$ and
$\tilde{z}_{k,n}=x_{k,n}\sqrt{\frac{q}{np}}$. Then, by arguments similar
to (\ref{expanA_1}) and (\ref{expanA_2}),
\begin{eqnarray*}
A_3&=&\left(nq-x_{k,n}\sqrt{npq}+\frac{1}{2}\right)\log
\left(1-x_{k,n}\sqrt{\frac{p}{nq}}\right)\nonumber\\
	&=& -\tilde{y}_{k,n}\left(nq+\frac{1}{2}\right)
	+\frac{\tilde{y}^2_{k,n}}{2}\left(nq-\frac{1}{2}\right)+r_{5n}(k)
\end{eqnarray*}
and
\begin{eqnarray*}
A_4&=&\left(np+x_{k,n}\sqrt{npq}+\frac{1}{2}\right)\log
\left(1+x_{k,n}\sqrt{\frac{q}{np}}\right)\nonumber\\
	&=&\tilde{z}_{k,n}\left(np+\frac{1}{2}\right)
	+\frac{\tilde{z}^2_{k,n}}{2}\left(np-\frac{1}{2}\right)+r_{6n}(k)
\end{eqnarray*}
where for all $k$ and $n$ satisfying $|\tilde{y}_{k,n}|\vee
|\tilde{z}_{k,n}| \leq \delta$, 
\begin{eqnarray}
|r_{5n}(k)|+|r_{6n}(k)|& \leq &
\frac{n}{2}\left[q{|\tilde{y}_{k,n}|}^3+p{|\tilde{z}_{k,n}|}^3\right]
	+\frac{2}{{(1-\delta)}^3}\left[\left(nq+\frac{1}{2}
	+nq{|\tilde{y}_{k,n}|}\right){|\tilde{y}_{k,n}|}^3\right.\nonumber\\
	& & \qquad\left.+\left(np+\frac{1}{2}+np{|\tilde{z}_{k,n}|}
				\right){|\tilde{z}_{k,n}|}^3\right] .
\label{errbound3}
\end{eqnarray}
Hence, as in (\ref{logexpann}), it follows that 
\begin{eqnarray}
\log \left\{{n\choose k}p^kq^{n-k}\right\}&=& \epsilon^{**}-\log
\sqrt{2\pi npq}-\frac{1}{2}x^2_{k,n}+r_{7n}(k)
\label{bierrexp}
\end{eqnarray}
where for all $n,k$ satisfying 
$|\tilde{y}_{k,n}|\vee |\tilde{z}_{k,n}|
\leq \delta$,
\begin{eqnarray*}
|r_{7n}(k)|\leq
\bigg|\frac{1}{2}\left(\tilde{z}_{k,n}
+\tilde{y}_{k,n}\right)-\frac{1}{4}\left(\tilde{y}^2_{k,n}
+\tilde{z}^2_{k,n}\right)\bigg|+|r_{5n}(k)|+|r_{6n}(k)| .
\end{eqnarray*}
Note that
\begin{eqnarray*}
fq+fp+(1-f)p+(1-f)q&=&1,\\
{(fq)}^2+{(fp)}^2+{((1-f)p)}^2+{((1-f)q)}^2
&=&(1-2pq)(1-2(1-f))<1,
\end{eqnarray*}
 and by (3.4), 
$y_{k,n}=fqa_{k,n}$, $z_{k,n}=fpa_{k,n}$, $\tilde{y}_{k,n}=pa_{k,n}$,
and 
$\tilde{z}_{k,n}=qa_{k,n}$.
Hence, it follows that
\begin{eqnarray}
\frac{1}{2}\left(|y_{k,n}|+|\tilde{y}_{k,n}|+|z_{k,n}|
+|\tilde{z}_{k,n}|\right)+\frac{1}{4}\left(y^2_{k,n}
+\tilde{y}^2_{k,n}+z^2_{k,n}+\tilde{z}^2_{k,n}\right)\leq
\frac{1}{2}|a_{k,n}|+\frac{1}{4}a^2_{k,n} .
\label{midineq}
\end{eqnarray}
Now, combining (\ref{pmfbreak}), (\ref{logexpann}) and (\ref{errbound3})
and using (\ref{midineq})
 and the above identities,  after some algebra,
we get 
$$
\log P(k;n,M,N)=-\frac{x^2_{k,n}}{2(1-f)}-\frac{1}{2}\log
(2\pi npq(1-f))+ r^*_n(k),$$ where for all $k,n$ satisfying
$|a_{k,n}|\leq \delta$,
\begin{eqnarray}
|r^*_n(k)-\epsilon^*-\epsilon^{**}|
& \leq &
|r_{4n}(k)|+|r_{7n}(k)|\nonumber\\	
   &\leq &
\frac{npq}{2}{|a_{k,n}|}^3\left[(1-f){(fq)}^2+(1-f){(fp)}^2
+p^2+q^2\right]\nonumber\\
	& &
{}+\frac{2npq}{{(1-\delta)}^3}{|a_{k,n}|}^3\bigg[(1-f)f^2\left\{(1+\delta
fq)q^2+(1+\delta fp)p^2\right\}\nonumber\\
	& & {}+(1+\delta p)p^2+(1+\delta q)q^2\bigg]\nonumber\\
	& &
{}+\frac{2}{{(1-\delta)}^3}{|a_{k,n}|}^3\frac{1}{2}
\left[(1+f^3)(p^3+q^3)\right]+\frac{1}{2}|a_{k,n}|+\frac{1}{4}a^2_{k,n}
\nonumber\\	
&\leq  & 
\frac{1}{2}|a_{k,n}|+a^2_{k,n}\left\{\frac{1}{4}
+\frac{2\delta}{{(1-\delta)}^3}\right\}
+{|a_{k,n}|}^3npq\left(\frac{f}{4}+1\right)\left\{\frac{1}{2}
+\frac{2(1+\delta)}{{(1-\delta)}^3}\right\} .\nonumber\\
\label{longineq}
\end{eqnarray} 
Note that for all $k,n$  satisfying $|a_{k,n}|\leq \delta$, 
$$Np-k\geq
Np-(np+\delta (1-f)npq)> np\frac{(1-f)}{2}> 0$$
 and
$$
Nq-(n-k)> nq\frac{(1-f)}{2}> 0 .$$
Hence, by the error bound in Stirling's approximation, for all 
$k,n$ with $|a_{k,n}|\leq \delta$ and $6(np\wedge nq)\geq 1$,
\begin{eqnarray*}
\epsilon^{*}
& \geq &
\frac{1}{12Np+1}-\frac{1}{12(Np-k)}+\frac{1}{12Nq+1}
-\frac{1}{12(Nq-(n-k))}+\frac{1}{12(N-n)+1}-\frac{1}{12N}\nonumber\\
&\geq &
-\frac{12k+1}{(12Np+1)(12(Np-k))}-\frac{12(n-k)
+1}{(12Nq+1)(12(Nq-n+k))}\nonumber\\
	& \geq &
-\frac{1}{6Np(1-\delta)(1-f)}-\frac{1}{6Nq(1-\delta)(1-f)}\nonumber\\
	&=& -\frac{f}{6npq(1-\delta)(1-f)};
\end{eqnarray*}
\begin{eqnarray*}
\epsilon^* \leq 0+0+\left[\frac{1}{12(N-n)+1}
-\frac{1}{12N}\right]\leq
\frac{f}{6npq(1-\delta)(1-f)};
\end{eqnarray*}
\begin{eqnarray}
\epsilon^{**}& \leq &
\frac{1}{12n}-\frac{1}{12k+1}-\frac{1}{12(n-k)+1}
\leq 0 ;
\label{errbound4}
\end{eqnarray}
\begin{eqnarray*}
\epsilon^{**}&\geq &\frac{1}{12n+1}-\frac{1}{12k}-\frac{1}{12(n-k)}\geq
-\frac{n}{12k(n-k)}\geq -\frac{1}{6npq(1-\delta)}.
\end{eqnarray*}
Hence, the lemma follows from (\ref{longineq}) and the above
inequalities.\\[.2in]

\noindent
{\bf Lemma 3.2}~
Let $g:\mathbf{R}\longrightarrow [0,\infty)$ be 
 such that $g$
is $\uparrow$ on $(-\infty,a)$ and $g$ is $\downarrow$ on $(a,\infty)$
for some $a \in \mathbf{R}$. Then, for any $k\in \mathbf{N}$, $b \in
\mathbf{R}$ and $h \in (0,\infty)$,
\begin{eqnarray}
\sum_{i=o}^k g(b+ih) \leq \int_{b}^{b+hk} g(x)dx +2hg(x_0),
\label{lem2con}
\end{eqnarray}
where $g(x_0)=max\{g(b+ih):i=0,1,\ldots,k\}$.
\\[.2in]
{\bf{Proof}}: For $b\geq a$, by monotonicity, 
$$
h\sum_{i=0}^k
g(b+ih)\leq hg(b)+\int_b^{b+hk}g(x)dx \, .
$$
For $b<a$, let $k_1=\sup \{i:b+ih<a\}$ and $b_1=b+k_1h$. 
Then,
\begin{eqnarray*}
h\sum_{i=0}^{k_1} g(b+ih)&\leq&
\sum_{i=0}^{k_1-1}\int_{b+ih}^{b+(i+1)h}g(x)dx+hg(b+k_1h)\\
	&\leq & \int_{b}^{b_1}g(x)dx+hg(b_1).
\end{eqnarray*}
Hence, for $b<a$ and $k>k_1$, 
\begin{eqnarray*}
h\sum_{i=0}^{k}g(b+ih)&=&h\sum_{i=0}^{k_1}g(b+ih)
+h\sum_{i=k_1+1}^{k}g(b+ih)\\
	&=& h\sum_{i=0}^{k_1}g(b+ih)+h\sum_{j=0}^{k-k_1-1}g(b_1+h+jh)\\
	&\leq &\int_{b}^{b_1}g(x)dx
+hg(b_1)+hg(b_1+h)+\int_{b_1+h}^{b_1+h+(k-k_1-1)h}g(x)dx-hg(b_1)\\
	&\leq & \int_{b}^{b+hk}g(x)dx+2hg(x_0).
\end{eqnarray*}
For $b<a$ and $k<k_1$, it is easy to check (using the arguments above)
that bound (\ref{lem2con}) trivially holds. This completes the proof of
the lemma.
\\[.2in]

\noindent
{\bf Lemma 3.3}~
Let $\phi(x)=\frac{1}{\sqrt{2\pi}}\exp(-\frac{x^2}{2})$, $x \in
\mathbf{R}$. Then, for any $h \in (0,\infty)$, $b \in [0,\infty)$,  $j_0
\in \mathbf{N}$,
\begin{eqnarray}
&&\bigg|h\sum_{i=0}^{j_0}\phi(b+ih)
-\int_{b-\frac{h}{2}}^{b+(j_0+\frac{1}{2})h}\phi(x)dx\bigg|
\label{lem3con}\\
&\leq
&
\frac{h^2}{12}\left[\int_{b-\frac{h}{2}}^{b+j_0h
+\frac{h}{2}}|{\phi}^{''}(x)|dx
{}+(4+h)\max\left\{|{\phi}^{''}(x)|:b-\frac{h}{2}<x<b+j_0h
+\frac{h}{2}\right\}\right] .\nonumber
\end{eqnarray}

\vspace*{.1in}
\noindent
{\bf{Proof :}}~
Note that the function $|{\phi}{''}(x)|=|x^2-1|\phi(x)$ is
even, and on $[0,\infty)$, it is increasing on $[1,{3}^{1/2}]$ 
and decreasing on each of the intervals $[0,1)$ and 
$({3}^{1/2},
\infty)$, with the  maximum value
$\frac{1}{\sqrt{2\pi}}$ at $x=0$ and the   minimum value $0$ at $x=1$.
First suppose that $(b-\frac{h}{2}, b+(j_0+\frac 12)h)\cap
\{0,\sqrt{3}\} =\nul$. Then, 
writing $b_i=b+ih$, $i\geq 0$, 
and using Taylor's expansion, one can show that  
 the leftside of (\ref{lem3con}) is
bounded above by
\begin{eqnarray*}
\sum_{i=0}^{j_0}
\bigg|\int_{b_i-\frac{h}{2}}^{b_i+\frac{h}{2}}
\Big(\phi(x)-\phi(b_i)\Big)dx\bigg|&\leq
&
\frac{1}{2}\sum_{i=0}^{j_0}\int_{b_i-\frac{h}{2}}^{b_i
+\frac{h}{2}}{(x-b_i)}^2\left\{\sup_{y
\in (b_i-\frac{h}{2},b_i+\frac{h}{2})}|{\phi}^{''}(y)|\right\}dx\\
	&\leq&
\frac{1}{2}\sum_{i=0}^{j_0}\left(2\int_{0}^{\frac{h}{2}}y^2dy\right)
\times \left\{\bigg|{\phi}^{''}\left(b_i-\frac{h}{2}\right)\bigg|\vee
\bigg|{\phi}^{''}\left(b_i+\frac{h}{2}\right)\bigg|\right\}\\
	&\leq & \frac{h^3}{24}
\sum_{i=0}^{j_0}\left\{\bigg|{\phi}^{''}
\left(b_i-\frac{h}{2}\right)\bigg|+
\bigg|{\phi}^{''}\left(b_i+\frac{h}{2}\right)\bigg|\right\}\\
	&\leq & \frac{h^3}{12}\sum_{i=0}^{j_0+1}
\bigg|{\phi}^{''}\left(b_i-\frac{h}{2}\right)\bigg|\,\,.
\end{eqnarray*}
 Hence
by two applications of Lemma 3.2, one can show that
\begin{eqnarray*}
h\sum_{i=0}^{j_0+1}\bigg|\phi^{''}\left(b_i
-\frac{h}{2}\right)\bigg| &\leq &
\int_{b-\frac{h}{2}}^{b+j_0h+\frac{h}{2}}|\phi^{''}(x)|dx
+4\max\{|\phi^{''}(x)|:b-\frac{h}{2}\leq
x \leq b+j_0h+\frac{h}{2}\}.
\end{eqnarray*}
Next consider the case where $0\in [b-\frac{h}{2},
b+\frac{h}{2})$. Then, by Taylor's expansion,
$$
\left | h\phi(b) -\int_{b-\frac{h}{2}}^{b-\frac{h}{2}}
\phi(x)dx\right|
\leq h^3|\phi''(0)|/24.
$$
Now using  similar arguments for
the case `$\sqrt{3}\in(b-\frac{h}{2}, b+(j_0+\frac 12)h)
\neq \nul$' and using  
the above bounds,   
 one can complete the proof of the lemma.\\[.2in]
%
%
%
%
{\bf{Proof of Theorem 2.1:}} Suppose that (\ref{rcondn})
holds. Fix $\epsilon$\ $\in (0,1)$. By Chebyshev's inequality, for 
all $r \in \mathbf{N}$,
\begin{eqnarray}
P\left(\bigg|\frac{X_r-n_rp_r}{\sigma_r}\bigg|
>\frac{2}{\epsilon}\right)&\leq
& \frac{\epsilon^2}{4}\, .
\label{thm1eq1}
\end{eqnarray}
By Lemmas 3.1 and 3.3, for any  $r\in \IN$ with $f_r\leq
\frac{1}{2}$,
\begin{eqnarray*}
\Delta_{1r}(\epsilon)&\equiv
&\sup_{-\frac{2}{\epsilon}\leq a <b\leq
\frac{2}{\epsilon}} \bigg|P\left(a<\frac{X_r-n_rp_r}{\sigma_r}\leq
b\right)-\left[\Phi(b)-\Phi(a)\right]\bigg|\nonumber\\
	&\leq & \sum_{-\frac{2\sigma_r}{\epsilon}<k-n_rp_r\leq
\frac{2\sigma_r}{\epsilon}}\bigg|P(k;n_r,M_r,N_r)
-\frac{1}{\sigma_r}\phi\left(\frac{k-n_rp_r}{\sigma_r}\right)\bigg|\\
	& &\quad {}+\sum_{-\frac{2}{\epsilon}\leq a<b\leq
\frac{2}{\epsilon}}\bigg|\sum_{a\sigma_r<k-n_rp_r\leq
b\sigma_r}\frac{1}{\sigma_r}\phi\left(\frac{k-n_rp_r}{\sigma_r}\right)
-[\Phi(b)-\Phi(a)]\bigg|\\
	& \leq & \frac{C}{\sigma^2_r}\sum_{-\frac{2\sigma_r}{\epsilon}<
k-n_rp_r\leq
\frac{2\sigma_r}{\epsilon}}\exp\left(\frac{C}{\sigma_r}\right)
\exp\left(-\frac{{(k-n_rp_r)}^2}{\sigma^2_r}
\left[\frac{1}{2}-\frac{C}{\sigma_r}\right]\right)\\
	& &\qquad
{}+\frac{C}{\sigma^2_r}\left[\int_{-\infty}^{\infty}|\phi^{''}(x)|dx
+1\right]+\frac{2}{\sqrt{2\pi}\sigma_r}\\
	&\leq &
\frac{C}{\sigma_r}\left[\int_{-\infty}^{\infty}\exp\left(
-\frac{x^2}{4}\right)dx+1\right],
\end{eqnarray*}
provided $\frac{C}{\sigma_r}< \frac{1}{4}$. Hence, there exists an
$r_0\in \mathbf{N}$ such that for all $r\geq r_0$ with $f_r\leq
\frac{1}{2}$
 $$
\Delta_{1r}(\epsilon)<\frac{\epsilon}{4}.
$$
Also by Mill's ratio,
$\Phi(-\frac{2}{\epsilon})+1-\Phi(\frac{2}{\epsilon})<\epsilon
\phi(\frac{2}{\epsilon})$. Hence, using (\ref{thm1eq1}) and the above
inequalities, it can be shown that for all $r\geq r_0$ with $f_r\leq
\frac{1}{2}$,
\begin{eqnarray}
\Delta_r(\epsilon)<\epsilon.
\label{thm2eq2}
\end{eqnarray}

Next suppose that $f_r>\frac{1}{2}$.
Consider the collection of $N_r-n_r$
objects that are left after the sample of size $n_r$ has been
selected from the population of size $N_r$. 
Let $Y_r=$the number of `type A'-objects in this 
collection.
Then, for all $r\in\IN$ and $j\in\mathbf{Z}$,
\begin{eqnarray}
Y_r\sim Hyp(N_r-n_r;M_r,N_r),\quad \mbox{and}
\quad P(X_r=j)=P(Y_r=M_r-j).
\label{newvardist}
\end{eqnarray}
Hence, 
$$
P(X_r \leq k)=\sum_{j=0}^k P(X_r=j)=\sum_{j=0}^k
P(Y_r=M_r-j)=P(Y_r\geq M_r-k).
$$
Further, note that $Var(Y_r)=(N_r-n_r)p_rq_r\left(1
-\frac{N_r-n_r}{N_r}\right)=\sigma^2_r$.
Hence, for each  $x\in \mathbf{R}$,
\begin{eqnarray*}
P\left(\frac{X_r-n_rp_r}{\sigma_r}\leq x\right)&=&P\left(X_r\leq
n_rp_r+x\sigma_r\right)\\
	&=& P(X_r\leq \lfloor n_rp_r+x\sigma_r \rfloor)\\
	&=& P(Y_r\geq M_r-\lfloor n_rp_r+x\sigma_r \rfloor)\\
	&=& P\left(\frac{Y_r-(N_r-n_r)p_r}{\sigma_r}\geq
\frac{M_r-\lfloor n_rp_r+x\sigma_r \rfloor
-(N_r-n_r)p_r}{\sigma_r}\right)\\
	&=& P(\tilde{Y}_r\geq \check{x}_r)\qquad\qquad \mbox{(say)},
\end{eqnarray*}
where $\tilde{Y}_r=\frac{Y_r-(N_r-n_r)p_r}{\sigma_r}$ and
 $\check{x}_r=\frac{M_r-\lfloor n_rp_r+x\sigma_r \rfloor
-(N_r-n_r)p_r}{\sigma_r}$. Note that,
$$\check{x}_r<\frac{1}{\sigma_r}\left[N_rp_r-(n_rp_r+x\sigma_r-1)
-N_rp_r+n_rp_r\right]=-x+\sigma_r^{-1}$$
and similarly, $\check{x}_r\geq -x$.
Hence, this implies, 
$$
P(\tilde{Y}_r<\check{x}_r)\leq P(\tilde{Y}_r\leq
\check{x}_r)\leq P(\tilde{Y}_r\leq -x+\sigma^{-1}_r)
$$ 
and 
$$
P(\tilde{Y}_r <
\check{x}_r)\geq P(\tilde{Y}_r<-x)\geq P(\tilde{Y}_r\leq
-x-\sigma^{-1}_r).
$$
Now using the above identity and  inequalities, we have
\begin{eqnarray}
\lefteqn{\bigg|P\left(\frac{X_r-n_rp_r}{\sigma_r}\leq
x\right)-\Phi(x)\bigg|=|P(\tilde{Y}_r\geq
\check{x}_r)-(1-\Phi(-x))|
=|\Phi(-x)-P(\tilde{Y}_r< \check{x}_r)|}\nonumber\\
      &\leq &  
	\max \{|P(\tilde{Y}_r\leq-x-\sigma^{-1}_r)
	-\Phi(-x-\sigma^{-1}_r)|,|P(\tilde{Y}_r\leq
-x+\sigma^{-1}_r)-\Phi(-x+\sigma^{-1}_r)|\}\nonumber\\
	& & \quad {}+\max\{|\Phi(-x)
	-\Phi(-x-\sigma^{-1}_r)|,|\Phi(-x)-\Phi(-x+\sigma^{-1}_r)|\}. 
\label{pineq1}
\end{eqnarray}
By repeating the arguments leading to (\ref{thm2eq2}), it follows that
there exists $r_1 \in \mathbf{N}$ such that $\forall r\geq r_1$ with
$(1-f_r)\leq \frac{1}{2}$,
\begin{eqnarray}
\sup_{x \in \mathbf{R}}|P(\tilde{Y}_r\leq x)-\Phi(x)|
\leq \epsilon .
\label{supineq}
\end{eqnarray}
Hence, (\ref{deltadef}) now follows from
(\ref{rcondn}),(\ref{thm2eq2}),(\ref{pineq1}) and (\ref{supineq}), with
$W\sim N(0,1)$. In particular, 
if (\ref{rcondn}) holds, then
one must have $\mu=0$ and
$\sigma=1$.

Conversely, suppose that (\ref{deltadef}) holds for some
$\mu \in \mathbf{R}$ and $\sigma \in (0,\infty)$. Then, for any
sequences $\{a_r\}_{r\geq 1}$,$\{b_r\}_{r\geq 1}$ $\subset
\mathbf{R}$ with $a_r<b_r$ for all $r\geq 1$,
\begin{eqnarray}
\bigg|P\left(a_r<\frac{X_r-n_rp_r}{\sigma_r}\leq b_r\right)-P(a_r<W\leq
b_r)\bigg|\leq 2\Delta_r \rightarrow 0\quad 
\mbox{as}\quad r \rightarrow
\infty.
\label{pineq2}
\end{eqnarray}
If possible, suppose that 
$\sigma_r<1$ infinitely often. Then,
 we can pick $a_r,b_r \in [-1,1]$
such that for all such $r$, $a_r-b_r=1$ and 
$$
\frac{\lfloor
n_rp_r\rfloor -n_rp_r}{\sigma_r}<a_r<b_r<\frac{\lfloor n_rp_r\rfloor 
+1-n_rp_r}{\sigma_r}.
$$
Then,
$$
P\left(a_r<\frac{X_r-n_rp_r}{\sigma_r}\leq b_r\right)=0
$$
 but
$$
P\left(a_r<W\leq b_r\right)\geq \inf\{P(a<W\leq b):a,b \in
[-1,1],b-a=1\}>0,
$$ infinitely often. This  contradicts 
(\ref{pineq2}). Hence, we may suppose that $\sigma_r\geq 1$ 
 for all but 
finitely many $r$'s.

Now define
$a_r=\frac{\lfloor n_rp_r\rfloor -n_rp_r+\frac{1}{3}}{\sigma_r}$ and
$b_r=\frac{\lfloor n_rp_r\rfloor -n_rp_r+\frac{2}{3}}{\sigma_r}$. 
Since $P(X_r\in \{0,1,\ldots,n_r\})=1$,
\begin{eqnarray*}
P\left(a_r<\frac{X_r-n_rp_r}{\sigma_r}\leq
b_r\right)=P\left(\lfloor n_rp_r\rfloor+\frac{1}{3}<X_r\leq \lfloor
 n_rp_r\rfloor+\frac{2}{3}\right)=0.
\end{eqnarray*}
Next using the definitions of $a_r$, $b_r$, and the fact that
`$x-1<\lfloor x\rfloor \leq x$ for all $x\in\mathbf{R}$', 
we get
\begin{eqnarray}
-\frac{2}{3\sigma_r}<a_r<b_r\leq \frac{2}{3\sigma_r},\quad r\geq 1.
\label{pineq3}
\end{eqnarray}
By (\ref{pineq2}) and (\ref{pineq3}), it follows that 
\begin{eqnarray*}
\frac{1}{3\sigma_r}\min \{\phi_{\sigma}(x-\mu):|x|\leq
\frac{2}{3\sigma_r}\}&\leq & \int_{a_r}^{b_r}\phi_{\sigma}(x-\mu)dx\\
	&=& P(a_r<W\leq b_r)\\
	&=& \bigg|P\left(a_r<\frac{X_r-n_rp_r}{\sigma_r}\leq
b_r\right)-P(a_r<W\leq b_r)\bigg|\\
	&\longrightarrow & 0\quad \mbox{as}\quad r \rightarrow \infty.
\end{eqnarray*}
As a result, $\sigma_r\raw \infty$ as $r \raw \infty$
and (2.5) holds. This completes the proof of the theorem.\\[.2in]

{\it To ensure  economy of space, we shall first give a proof of 
Theorem 2.3 and then outline the main steps in the 
proof of Theorem 2.2.}\\[.2in]

\noindent
{\bf Proof of Theorem 2.3:}~ Let $r\in \mathbf{N}$ be an
integer such that (2.13) holds. Since $r$ will
be held {\it fixed} all through the proof, we shall drop 
$r$ from the notation for simplicity, and write 
$f_r =f$, $\sigma_r=\sigma$, $p_r = p$, $q_r = q$, $n_r-n$, etc.
First, suppose that $f\leq \frac{1}{2}$. Consider the
case $x\leq 0$. Let
$\tilde{x}_{k}=\frac{x_k}{\sqrt{1-f}}=\frac{k-np}{\sigma}$,
$k=0,1,\ldots,n$.
Define 
\begin{eqnarray*}
K_0&=&\sup\{k\in \mathbf{Z}_+:\tilde{x}_k\leq 0\}\\
K_1&=&\sup\{k\in \mathbf{Z}_+:\tilde{x}_k\geq -1\}\\ 
K_2&=&\sup\{k\in \mathbf{Z}_+:\tilde{x}_k\geq 
-\delta\sigma\}\quad \mbox{and}\\
J_x&=&\lfloor np+x\sigma \rfloor, x\in \mathbf{R},\\
\end{eqnarray*}
where $\delta \equiv \delta_r\in (0,\frac{1}{2}]$
is as in (2.12). 
Note that by definition,
\begin{eqnarray*}
K_1-1<np-\sigma\leq K_1,\quad{}&&
K_2-1<np-\delta\sigma^2 \leq K_2,\\
 \tilde{x}_j \in [-1,0]
\forall K_1\leq j \leq K_0 && \mbox{and}\quad
\tilde{x}_j \in
[-\delta\sigma,-1)\forall K_2\leq j< K_1.
\end{eqnarray*}
Hence, for any $x\in [-\delta\sigma,0]$, 
\begin{eqnarray}
&&\bigg|P\left(\frac{X-np}{\sigma} \leq x \right)-\Phi(x)\bigg|
=
|P(X\leq J_x)-\Phi(x)|\nonumber\\
	&\leq &
P(X<K_2)+\sum_{j=K_2}^{J_x}\bigg|P(X=j)
-\frac{\phi(\tilde{x}_j)}{\sigma}\bigg|
+\bigg|\sum_{j=K_2}^{J_x}\frac{\phi(\tilde{x}_j)}{\sigma}
-\Phi(x)\bigg|\nonumber\\
	&=& I_1+I_2+I_3,\qquad \mbox{say.}
\label{pineq4}
\end{eqnarray}
Consider $I_2$ for $x \in [-\delta\sigma,-1)$. Note that for $x<-1$,
$\frac{J_x-np}{\sigma} \leq x <-1$. Hence $J_x<K_1$ and 
$\tilde{x}_j<-1$ for all $j<J_x$.
From Lemma 3.1,
\begin{eqnarray}
|r^*(j)|& \leq & \frac{1}{6\sigma^2
(1-\delta)}+\left[\frac{{|\tilde{x}_j|}^2}{2\sigma}
+\frac{{|\tilde{x}_j|}^2}{\sigma^2}\left\{\frac{1}{4}+\frac{2\delta}{{(1
-\delta)}^3}\right\}
+\frac{{|\tilde{x}_j|}^3}{2\sigma}A\right]\nonumber\\
	&\equiv & r^{**}(j),
\label{pineq5}
\end{eqnarray} 
where $A=a_1\left(1+\frac{4(1+\delta)}{{(1-\delta)}^3}\right)$ and
$a_1\equiv a_{1r}=\frac{f+4}{4(1-f)}$ (cf. (2.12)).
%
 For  the given 
choice of $\delta$, it is easy to verify that
$\delta \leq \frac{1}{20}$ and 
 $\delta A<.59$. Hence
\begin{eqnarray}
|r^*(j)|& \leq & (0.2)
\sigma^{-2}+\frac{\tilde{x}^2_j}{2}\left[\frac{1}{\sigma}
+\frac{2}{\sigma^2}(0.3667)+\delta
A\right]\nonumber\\
	&\leq & (0.2) 
\sigma^{-2}+\frac{\tilde{x}^2_j}{2} \left[\min
\{0.86,\frac{6}{5\sigma}+0.59\}\right] .
\label{pineq7}
\end{eqnarray}
%
%
Now, from (\ref{pineq5}), for all $K_2\leq j< K_1$,
\begin{eqnarray}
|r^*(j)|& \leq & (0.2)
\sigma^{-2}+{|\tilde{x}_j|}^3\bigg|\left[\frac{1}{2\sigma}
+\frac{1}{\sigma^2}(0.3667)+\frac{3a_1}{\sigma}\right]
	\leq  4{|\tilde{x}_j|}^3\frac{a_1}{\sigma}.
\label{pineq8}
\end{eqnarray}

Next note that 
$$
\frac{J_x-np}{\sigma}\leq x \in \mathbf{R},
$$
$$
\int_{a}^{\infty}
y^3\exp(-\frac{by^2}{2})dy=\frac{1}{2b^2}(1+ba^2)e^{-ba^2}
\forall a,b\in (0,\infty),
$$ and that for any $a \in (0,\infty)$,
 the function
$g(y;a)=y^3\exp(-ay)$, $y \in [0,\infty)$, is increasing on
$[0,\sqrt{\frac{3}{2a}}]$, and decreasing on
$(\sqrt{\frac{3}{2a}},\infty)$. Hence, 
%
%
 by Lemmas 3.1 and 3.2, (\ref{pineq7}) and
(\ref{pineq8}), with $c=.07$, we have 
\begin{eqnarray}
I_2 &\leq &
\sum_{j=K_2}^{J_x}\bigg|\frac{\phi(\tilde{x}_j)}{\sigma}\exp(r^*(j))
-\frac{\phi(\tilde{x}_j)}{\sigma}\bigg|\nonumber\\	
&\leq&
\frac{1}{\sigma}\sum_{j=K_2}^{J_x}
\phi(\tilde{x}_j)|r^*(j)|\exp(|r^*(j)|)\nonumber\\
	&\leq &
\frac{4a_1}{\sqrt{2\pi}\sigma^2}\exp(\sigma^{-2})
\sum_{j=K_2}^{J_x}{|\tilde{x}_j|}^3
	\exp(-c\tilde{x}^2_j)\nonumber\\
& \leq &
\frac{4a_1\exp(\sigma^{-2})}{\sqrt{2\pi}\sigma}
\Big[\int_{\frac{K_2-np}{\sigma}}^{\frac{J_x
-np}{\sigma}}{|y|}^3\exp(-c|y|)dy\nonumber\\
&&\quad{}+\frac{2}{\sigma}\max\{{|y|}^3\exp(-c|y|)
: K_2\leq np+\sigma y \leq J_x\}\Big] \nonumber\\
%
&\leq& \frac{C}{\sigma(1-f)}\left[(1+x^2)\exp(-cx^2)\right] .
\label{pineq9}
\end{eqnarray}
Also, for $-1\leq x \leq 0$, by Lemma 3.1, 
\begin{eqnarray*}
\Delta_1(x)
&\equiv &
	\bigg| P\left(-1\leq \frac{X-np}{\sigma}\leq x\right)
	-\sum_{j=K_1}^{K_0}\frac{1}{\sigma}\phi(\tilde{x}_j)\bigg|\\
&\leq &
	 \sum_{j=K_1}^{K_0}
\bigg| P(X=j)-\frac{1}{\sigma}\phi(\tilde{x}_j)\bigg|\\
	&\leq &
\sum_{j=K_1}^{K_0}\exp\left(
-\frac{\tilde{x}^2_j}{2}\right)|r^*(j)|
\frac{\exp(|r^*(j)|)}{\sqrt{2\pi}\sigma}.
\end{eqnarray*}
For $K_1\leq j\leq K_0$, from (\ref{pineq5}) and (\ref{pineq7}),
\begin{eqnarray*}
|r^*(j)|&\leq &
\left[\frac{1}{2\sigma}|\tilde{x}_j|+r^{**}(j)\right]\wedge
\left[\frac{1}{5\sigma^2}+\frac{1}{2\sigma}
+\frac{1}{2\sigma^2}(0.3667)+\frac{A}{2\sigma}\right]\\
	&\leq &
\left[\frac{1}{2\sigma}+\frac{1}{5\sigma^2}
+(0.43)\tilde{x}^2_j\right]\wedge
\left[\frac{1}{2\sigma}+\frac{1}{5\sigma^2}+\frac{0.3667}{\sigma^2}
+\frac{3a_1}{\sigma}\right]\\
	&\leq & \frac{1}{\sigma}+\left[(.43)\tilde{x}^2_j\right]\wedge
\left[\frac{4a_1}{\sigma}\right] .
\end{eqnarray*}
Hence, for $-1\leq x\leq 0$, noting that $K_0-K_1\leq \sigma$,
\begin{eqnarray}
|\Delta_1(x)|&\leq &
\sum_{j=K_1}^{K_0}\exp(-\tilde{x}^2_j(0.07))
\exp(\sigma^{-1})\frac{5a_1}{\sqrt{2\pi}\sigma^2}\nonumber\\
	&\leq &
(K_0-K_1)\exp(\sigma^{-1})\frac{5a_1}{\sqrt{2\pi}\sigma^2}\nonumber\\
	&\leq & \frac{C}{\sigma} .
\label{pineq10}
\end{eqnarray}
Thus, the bound (\ref{pineq9}) on $I_2$
holds for all $x \in [-\delta\sigma,0]$.

Next consider $I_1$. Note that for $j\in\{0,1,\ldots,n\}$, 
\begin{eqnarray}
\lefteqn{P(X=j+1)>=<P(X=j)}\hspace{1.8in}\nonumber\\
\Leftrightarrow \frac{Np-j}{j+1}.\frac{n-j}{Nq-n+j+1}>&=&<1 \nonumber\\
\Leftrightarrow j<=>np-\frac{Nq+1}{N+2}\, .
\label{pineq11}
\end{eqnarray}
Thus, $P(X=j)<P(X=j+1)$ for all $0\leq j \leq np-1$. Hence, by
(\ref{pineq7}) and Lemma 3.1,
\begin{eqnarray*}
I_1&=&\sum_{j=0}^{K_2-1}P(X=j)\\
	&< & K_2 P(X=K_2)\\
	&\leq &
K_2\frac{1}{\sigma}\phi\left(\tilde{x}_{K_2}\right)\exp(r^*(K_2))\\
	&\leq &
\frac{K_2}{\sqrt{2\pi}\sigma}
\exp\left(\frac{1}{5\sigma^2}\right)\exp(-\tilde{x}^2_{K_2}(.07))\\
	&\leq &
\frac{K_2}{\sqrt{2\pi}\sigma}\exp\left(\frac{1}{5\sigma^2}\right)
\exp\left(-{\left(\delta
\sigma-\frac{1}{\sigma}\right)}^2(0.07)\right)\\
	&\leq &
\frac{K_2}{\sqrt{2\pi}\sigma}\exp(-\delta^2\sigma^2(0.07)+2\delta
(0.07)+0.13\sigma^{-2})\\
	&\leq &
\frac{np}{\sqrt{2\pi}\sigma}\exp(-\delta^2\sigma^2(0.07))\exp(0.014)\\
	&\leq & (q(1-f))^{-1}\sigma \exp(-\delta^2\sigma^2(0.07)).
\end{eqnarray*}
It is easy to check that,
\begin{displaymath}
\frac{\sigma \exp(-\delta^2\sigma^2(0.07))}{(1+x^2)\exp(-x^2(0.07))}\leq
\left\{\begin{array}{c@{\quad:\quad}l} \frac{2}{(0.07)\delta^2\sigma} &
\mbox{if} \quad x \in [0,\frac{\delta\sigma}{\sqrt{2}}], \\
\frac{2}{\delta^2 \sigma} & 
\mbox{if} \quad x \in [\frac{\delta\sigma}{\sqrt{2}},\delta\sigma] .
\end{array}\right.
\end{displaymath}
Hence, it follows that for all $x \in [-\delta a,0]$,
\begin{eqnarray}
I_1\leq \frac{C}{\delta^2q\sigma(1-f)}(1+x^2)\exp(-x^2(0.07)).
\label{pineq12}
\end{eqnarray}
Next note that by definition, $\tilde{x}_{J_x}
\leq x$ and $\tilde{x}_{K_2}\leq -\delta\sigma
+\sigma^{-1}$. Hence, for
$x \in [-\delta\sigma,0]$, by Lemma 3.3, 
\begin{eqnarray*}
I_3&\leq &
\bigg|\frac{1}{\sigma}\sum_{j=K_2}^{J_x}\phi(\tilde{x}_j)
	-\int_{\tilde{x}_{K_2}-{(2\sigma)}^{-1}}^{\tilde{x}_{J_x}
	+{(2\sigma)}^{-1}}\phi(y)dy\bigg|+\bigg|\Phi(x)
	-\Phi\left(\tilde{x}_{J_x}+{(2\sigma)}^{-1}\right)\bigg|
	+\Phi\left(\tilde{x}_{K_2}-{(2\sigma)}^{-1}\right)\\
	&\leq &
\frac{1}{12\sigma^2}\left[\int_{-\infty}^{x+\frac{1}{2\sigma}}
|\phi^{''}(y)|dy+5\max \{|\phi^{''}(y)|:
-\infty<y<x+\frac{1}{2\sigma}\}\right]\\
&&\qquad {}
+\Phi\left(x+\frac{1}{2\sigma}\right)-
\Phi\left(x-\frac{1}{2\sigma}\right)
	+\Phi(-\delta\sigma+\frac{1}{2\sigma}).
\end{eqnarray*}
Note that for any $a\in
(0,\infty)$,
$$
\int_{a}^{\infty}y^2e^{\left(-\frac{y^2}{2}\right)}dy 
\leq \frac{1}{a}\int_{a}^{\infty}y^3e^{\left(-\frac{y^2}{2}\right)}dy
=\frac{2}{a}\int_{\frac{a^2}{2}}^{\infty}te^{-t}dt
=\frac{a^2+2}{a}e^{-\frac{a^2}{2}};
$$
$$
\int_{a}^{\infty}y^2e^{-\frac{y^2}{2}}dy\leq
\int_{0}^{\infty}y^2e^{-\frac{y^2}{2}}dy\leq \sqrt{\frac{\pi}{2}};
$$
$$
\max\{|\phi^{''}(y)|:a<y<\infty\}\leq
\frac{1}{\sqrt{2\pi}}I(0<a<\sqrt{3})
+|\phi^{''}(a)|I(a\geq \sqrt{3});
$$
$$\exp\left(-\frac{{(a-{(2\sigma)}^{-1})}^2}{2}\right)\leq
\exp\left(-\frac{a^2}{2}+\frac{a}{2\sigma}\right)\leq
\exp\left(-\frac{a^2}{2}+\frac{\delta}{2}\right),\forall a \in
(0,\delta\sigma).$$ 
Also note that, for $0<a\leq 1$, $b\in (0,\infty)$,
$$
1-\Phi(b)\leq \frac{1}{b}\phi(b),
$$
$$ 
1-\Phi(a)\leq \int_{a}^{1}\phi(x)dx+\phi(1)\leq \phi(a)(1-a)+\phi(a)
=(2-a)\phi(a).$$
Thus, for any $x\in (0,\infty)$, 
$$
\Phi(x)\leq e^{-\frac{x^2}{2}}.
$$
Since  ${(2\sigma)}^{-1}<\frac{1}{8}$ and
$|y+{(2\sigma)}^{-1}|\leq |y|$  for
$y <-\frac{1}{8}$, we have,  for all $x \in
[-\delta a,0]$,
\begin{eqnarray}
I_3
&\leq &
\frac{1}{12\sigma^2}\Big[2I(-2\leq x\leq 0)
+5|x|\phi(x+{(2\sigma)}^{-1})I(-\delta\sigma 
		\leq x \leq -2) \nonumber\\
& & 
{}+5\left\{\frac{1}{\sqrt{2\pi}}
I(-2\leq x\leq 0)+(x^2+1)\phi(x+{(2\sigma)}^{-1})
I(-\delta\sigma \leq x \leq-2)\right\}\Big]\nonumber\\
& &
 {}+\frac{1}{\sqrt{2\pi}\sigma}I(-2\leq x \leq
0)+\frac{1}{\sigma} \phi\left(x+{(2\sigma)}^{-1}\right)I(-\delta a\leq
x<-2)+\Phi(-\delta\sigma+{(2\sigma)}^{-1}) \nonumber\\
&\leq & 
\frac{1}{2\sigma}
I(-2\leq x\leq 0)+2 \left\{\frac{x^2+1}{2\sigma^2}
+\frac{1}{\sigma}\right\}\phi\left(x
+\frac{1}{2\sigma}\right)I(-\delta\sigma
\leq x \leq
-2)+\exp\left(-\frac{{(\delta\sigma
-\frac{1}{2\sigma})}^2}{2}\right)\nonumber\\
&\leq &
 \frac{C}{\sigma}(1+|x|)\exp\left(-\frac{x^2}{2}\right)
\, .
\label{pineq13}
\end{eqnarray}
Next note that
 $$
P\left(\frac{X-np}{\sigma}\leq x\right)=0\forall
x<-\frac{np}{\sigma}
$$ 
and for $-\frac{np}{\sigma}\leq x \leq
-\delta\sigma$, 
\begin{eqnarray*}
P\left(\frac{X-np}{\sigma}\leq x\right)\leq I_1 
&\leq &
{(q(1-f))}^{-1}\sigma\exp(-\delta^2\sigma^2(0..07))\\
	&=&
{(q(1-f))}^{-1}\frac{{(\delta\sigma)}^2}{\delta^2\sigma}
\exp\Big(-\delta^2q^2{(1-f)}^2
{\left[\frac{-np}{\sigma}\right]}^2(0.07)\Big)\\
	&\leq & 
{(\delta^2q(1-f)\sigma)}^{-1}{|x|}^2\exp\Big(-\delta^2q^2{(1
-f)}^2x^2(0.07)\Big).
\end{eqnarray*}
Hence, for all $x\leq -\delta\sigma$,
\begin{eqnarray}
\bigg|P\left(\frac{X-np}{\sigma}\leq x\right)-\Phi(x)\bigg|
&\leq &
\frac{{|x|}^2\exp(-\delta^2q^2{(1-f)}^2x^2(0.07))
+\exp\left(-\frac{x^2}{2}\right)}{\delta
q (1-f)\sigma}\nonumber\\
	&\leq & \frac{2}{\delta q
(1-f)\sigma}x^2\exp\Big(-\delta^2q^2{(1-f)}^2x^2(0.07)\Big).
\label{pineq14}
\end{eqnarray}
Now using the fact that $\delta \in
\left[\frac{1}{25},\frac{1}{20}\right]$  for all 
$f \in (0,\frac{1}{2}]$,
from (\ref{pineq9}),(\ref{pineq10}) and (\ref{pineq12})-(\ref{pineq14}),
it follows that there exist numerical constants $C_1$ and $C_2$, not
depending on $n,M,N$,  such that for all $x \in (-\infty,0]$,
$$
\bigg|P\left(\frac{X-np}{\sigma}\leq x\right)-\Phi(x)\bigg|\leq
\frac{C_1}{\sigma q}(1+x^2)\exp(-C_2qx^2),
$$
 provided $\delta\sigma>1$. This proves 
(\ref{main3bound}) for $x \in (-\infty,0]$ and $f\leq
\frac{1}{2}$.

To prove the theorem for $x\geq 0$ and $f\leq \frac{1}{2}$, 
define
$$
V_r=n_r-X_r,\quad r \in \mathbf{N}.
$$
Note that $V_r$ has a  Hypergeometric distribution with parameters
$n_r,N_r-M_r,N_r$. Further,
$$
\frac{X_r-n_rp_r}{\sigma_r}
=-\frac{V_r-n_rq_r}{\sigma_r}\forall r\in {\mathbf{N}}.
$$
Hence, the derived bound on the right tails of
$\frac{X_r-n_rp_r}{\sigma_r}$, can be obtained by repeating the
arguments above with $X_r$ replaced by $V_r$ and $p_r$ replaced by
$q_r$ for any $r$ such that $\delta\sigma_r> 1$. This
proves (\ref{main3bound}) for $x \in [0,\infty)$ and $f \leq
\frac{1}{2}$. The proof of (\ref{main3bound}) for `$f \in
[\frac{1}{2},1]$ and $x\in \mathbf{R}$'
follows by replacing the above arguments with $X_r,f_r$
replaced by $Y_r,1-f_r$ respectively and using the bound
(\ref{newvardist}) and (\ref{pineq1}). This completes 
the prrof of the
theorem.\\[.2in]

\noindent
{\bf{Proof of Theorem 2.2:}} As in the proof of Theorem 2.3,
 first we suppose that  $f_r\leq \frac{1}{2}$. 
By (3.1),  (3.32), (3.36), (3.37), and (3.40),
  it
follows that for all $r$ with $\delta_r\sigma_r> 1$, 
\begin{eqnarray}
\sup_{x\in [-\delta_r\sigma_r,0]}|\Delta_r^*(x)|
&\leq &
P(X_r<K_2)+\sup_{x \in [-\delta_r\sigma_r,0]}\{I_2+I_3\}\nonumber\\
	&\leq & P(X_r\leq K_2-1)+\frac{C}{\sigma_r}.
\label{pineq15}
\end{eqnarray}
By Chebyshev's inequality, noting that
$K_2-1<n_rp_r-\delta_r\sigma^2_r\leq K_2$, we have 
\begin{eqnarray}
P(X_r\leq K_2-1)
&\leq &
P\left(\bigg|\frac{X_r-n_rp_r}{\sigma_r}\bigg|\geq
\bigg|\frac{K_2-n_rp_r-1}{\sigma_r}\bigg|\right)\nonumber\\
&\leq & 
\frac{Var(X_r)}{{(K_2-1-n_rp_r)}^2}\nonumber\\
	&\leq &
\frac{N_r\sigma^2_r}{N_r-1}{(\delta_r\sigma^2_r)}^{-2}\nonumber\\
	&\leq & \frac{2}{\delta^2_r\sigma^2_r}.
\label{pineq16}
\end{eqnarray}
Also, 
\begin{eqnarray}
\sup_{-\infty\leq x \leq -\delta_r\sigma_r}|\Delta^*_r(x)|&\leq &
P(X_r\leq K_2-1)+\Phi(-\delta_r\sigma_r)\nonumber\\
	&\leq & \frac{C}{\delta^2_r\sigma^2_r}.
\label{pineq17}
\end{eqnarray}
Since $\delta_r\geq \frac{1}{22.5}$ for all $r$ with $f_r \leq
\frac{1}{2}$, from (\ref{pineq15})-(\ref{pineq17}),
it follows that there
exists a universal constant $C_3$ such that 
for all $r$ with $\delta_r\sigma_r> 1$ 
and $f_r\leq \frac{1}{2}$, 
$$
\sup_{x \leq 0}|\Delta^*_r(x)|\leq \frac{C_3}{\sigma_r}.
$$
Now retracing the arguments in the proof of Theorem 2.3 
for the case ``$x\geq 0,f_r\leq \frac{1}{2}$" 
(with the variable $V_r$) and for  the case
``$x\in \mathbf{R},f>\frac{1}{2}$" 
(with $Y_r$), one can complete the
proof of Theorem 2.3. \\
\\
{\bf{Proof of Corollary 2.4}}: Use (\ref{main3bound}) and the 
inequality ``$\exp(x)\geq (1+x)$  for all $x \in (0,\infty)$''.
\\[.2in]
\begin{center}
{\bf{\large REFERENCES}}
\end{center}
\begin{description}
\item
Babu, G.J. and Singh, K. (1985).  
Edgeworth expansions for sampling without replacement from 
finite populations. {\it  Journal of Multivariate Analysis}  
{\bf 17}   261-278.
\item
Bloznelis, M.   (1999). A Berry-Esseen bound for finite
 population Student's statistic. {\it Annals of Probability}
 {\bf  27}   2089-2108.
\item
Bloznelis and G\"{o}tze, F. (2000).
An Edgeworth expansion for finite-population $U$-statistics.
{\it Bernoulli}  {\bf 6}   729-760.
\item
Chen, J. and Sitter, R.R. (1993). Edgeworth expansion and 
the bootstrap for stratified sampling without 
replacement from a finite population. {\it 
 Canadian Journal of Statistics} {\bf  21}  347-357.
\item
Erdos, P. and Renyi, A. (1959).
On the Central Limit Theorem for samples from a
 finite population.  {\it Magyar Tudoanyos 
Akademia Budapest
Matematikai Kutato Intezet Koezlemenyei}, Trudy 
Publications, {\bf 4} 49-57.
\item
Feller, W. (1971).  {\it An introduction to probability theory
and its applications. Volume I}. Wiley, New York, NY.
\item
Hajek, J. (1960).  Limiting distributions in simple random 
sampling from a finite population.  {\it Magyar Tudoanyos 
Akademia Budapest
Matematikai Kutato Intezet Koezlemenyei}, Trudy 
Publications, {\bf 5} 361-374.
\item
Kokic, P. N.   and  Weber, N. C.  (1990). 
An Edgeworth expansion for $U$-statistics based on 
samples from finite populations. {\it Annals of Probability}
{\bf   18}  390-404.
\item
Nicholson,  W.L. (1956). On the Normal approximation to
 the Hyprgeometric distribution. {\it Annals of Mathematical
Statistics} {\bf 27} 471-483.
\item
Madow, W.G. (1948). On the limiting distributions
of estimates based on samples from finite universes.
{\it Annals of Mathematical
Statistics} {\bf 19} 535-545.
\end{description}
\end{document}